\documentclass{article}
\usepackage{amssymb}
\usepackage{amsmath}
\usepackage{amsthm}
 \newcommand{\nn}{\mathbb N} 
\newcommand{\zz}{\mathbb Z} \newcommand{\rr}{\mathbb R}

\newcommand{\spa}{\hbox{\rm span}}

\newcommand{\cl}{\hbox{\rm cl}}

\newtheorem{thm}{Theorem}
\newtheorem{prop}{Proposition}

\newtheorem{exam}{Example}
\newtheorem{coro}{Corollary}
\newtheorem{asser}{Assertion}
\newtheorem{definition}{Definition}

\begin{document}
\normalsize

\title{\Large\bf Archimedeanization of ordered vector spaces}
\author{\bf Eduard Yu. Emelyanov\footnote{Middle East Technical University, Ankara, Turkey}}
\maketitle

{\bf Abstract:} {\normalsize 
In the case of an ordered vector space (briefly, OVS) with an order unit, the Archimedeanization 
method has been developed recently by Paulsen and Tomforde \cite{PT}. We present a general version of the 
Archimedeanization which covers arbitrary OVS.}



\section{Preliminaries}

In the present paper, we always exclude 0 from natural numbers ($\nn:=\{1,2,...\}$)
and consider vector spaces only over reals.\\

A subset $W$ of a vector space $V$ is called a {\it wedge} if it satisfies 
$$
  W+W\subseteq W
$$ 
and
$$
  rW\subseteq W\ \ \ \ (\forall\ r\ge 0).
  \eqno(1)
$$ 
A wedge $K$ is called a {\it cone} if 
$$
  K\cap(-K)=\{0\}.
$$ 
A cone $K$ is said to be {\it generating} if 
$$
  K-K=V.
  \eqno(2)
$$ 
Given a cone (resp., a wedge) $V_+$ in $V$, we say that $(V,V_+)$ is an {\it ordered vector space} 
(briefly an OVS) (resp., a {\it pre-ordered vector space} or a POVS). The partial ordering (resp., partial pre-ordering)
$\le$ on $V$ is defined by 
$$ 
  x\le y \ \ \text{\rm whenever} \ \ y-x\in V_+\,.
$$ 
$(V,V_+)$ is also denoted by $(V,\le)$ or simply by $V$, if the partial ordering $\le$ is well understood.
We say that $e\in V$ is an {\it order unit} if for each $x\in V$ there is some $n\in\nn$ such that
$x\le ne$. Clearly, $V_+$ is generating if $V$ possesses an order unit.

For every $x,y$ in a POVS $(V,\le)$, the (possibly empty, if $x\not\le y$) set 
$$
  [x,y]:=\{z:x\le z\le y\}
$$ 
is called an {\it order interval}. A subset $A\subseteq V$ is said to be {\it order convex} 
if $[a,b]\subseteq A$ for all $a,b\in A$. Any order convex vector subspace of $V$ is called an {\it order ideal}. 
If $Y$ is an order ideal in a POVS $V$, the quotient space $V/Y$ is partially ordered by:
$$
  [0]\le [f] \ \ \text{\rm if} \ \ \exists q\in Y \ \text{\rm with} \ \  0\le f+q\ .
$$
The order convexity of $Y$ is needed for the property
$$
  [f]\le [g] \le [f] \ \Rightarrow [f]=[g]\,,
$$
which guarantees that $[V_+]_Y:=V_++Y$ is a cone in  $V/Y$.

Let $x_n\in V$ and $u\in V_+$. The sequence $(x_n)_n$ is said to be {\it $u$-uniformly convergent}
to a vector $x\in V$, in symbols $x_n\stackrel{(u)}{\to}x$, if 
$$
  k_n(x_n-x)\in[-u,u]
$$ 
for some sequence $k_n\uparrow\infty$ in $\nn$. Clearly, $x_n\stackrel{(u)}{\to}x$ iff  
$$
  x_n-x\in[-\varepsilon_nu,\varepsilon_nu]
$$ 
for some sequence $\varepsilon_n$ of reals such that $\varepsilon_n\downarrow 0$. A subset $A\subseteq V$ is said to be 
{\it uniformly closed} in $V$ if it contains all limits of all $u$-uniformly convergent sequences $(x_n)_n$ in $A$ for all $u\in V_+$.

\begin{definition}\label{D1}
A wedge $K$ in a vector space $V$ is said to be {\rm Archi\-me\-dean} $($also we say that the POVS $(V,K)$ is {\rm Archi\-medean}$)$, whenever
$$
  (\forall x,y\in V)\bigg[[(\forall n\in\nn) \ x-ny\in K]\Rightarrow [y\in -K]\bigg]\,.
  \eqno(3)
$$
\end{definition}
\noindent
Clearly, $(3)$ in Definition \ref{D1} can be replaced by 
$$
  (\forall x\in V_+,y\in V)\bigg[[(\forall n\in\nn) \ x-ny\in K]\Rightarrow [y\in -K]\bigg]\,.
  \eqno(3')
$$
It is obvious that any subspace of an Archimedean POVS is Archi\-medean. 
It is also worth to notice that if $V$ admits a linear topology $\tau$ for which 
a wedge $K\subseteq V$ is closed, then $K$ is Archi\-medean (cf. \cite[Lem.2.3]{AT}). 
In the case when $\tau$ is Hausdorff and a wedge $K\subseteq V$ has nonempty $\tau$-interior, 
$K$ is Archi\-medean if and only if it is $\tau$-closed (cf. \cite[Lem.2.4]{AT}). 

\begin{asser}\label{A1}
An OVS $(V,V_+)$ is Archimedean if and only if 
$$
  (\forall y\in V_+)\inf_{n\ge 1}\frac{1}{n}y=0\,.
$$
\end{asser}

Clearly, one may define the Archimedean property by saying that, whenever $x\in V_+$ 
and $y\in V$ with $0\le x+ny$ for all $n\in\nn$, then $0\le y$. This observation
leads to the following definition, that is motivated by the definition of an {\it Archimedean order unit} 
in \cite[Def. 2.7]{PT}.

\begin{definition}\label{D2}
If $(V,V_+)$ is a POVS and $x\in V_+$, we say that $x$ is an {\rm Archimedean element} if,
whenever $y\in V$ with $x+ny\in V_+$ for all $n\in\nn$, then $y\in V_+$.
\end{definition}
\noindent
Notice that $x=0$ is always an Archimedean element. Clearly, $V$ is Archimedean 
if and only if every element $x\in V_+$ is Archimedean. Furthermore, if $0\le z\le x$ and $x$ is Archimedean, then $z$ 
is Archi\-me\-dean. In particular, in any POVS $V$ with an Archimedean order unit, all elements of $V_+$ are Archimedean. 

\begin{definition}\label{D3}
A wedge $W$ in a vector space $V$ is said to be {\rm almost Archime\-dean} 
$($also we say that the POVS $(V,W)$ is {\rm almost Archi\-medean}$)$ if
$$
  (\forall x,y\in V)\bigg[[(\forall n\in\zz) \ x-ny\in W]\Rightarrow [y=0]\bigg]\,.
  \eqno(4)
$$
\end{definition}
\noindent
Clearly, any $x\in V$ satisfying $[(\forall n\in\zz) \ x-ny\in W]$ in (4) is indeed an element of $W$.
OVS which are not almost Archimedean given in the following example.

\begin{exam}\label{E0}
Let $V=\rr^2$ be ordered by the cones 
$$
  K_1:=\{(x,y)|x>0\}\cup\{0\}\,, \ \ K_2:=\{(x,y)|x>0\}\cup\{(0,y)|y\ge 0\}\,.
$$
The OVS $(V,K_2)$ is denoted by $(\rr^2,\le_{lex})$ and called the {\em Euclidean plane with
the lexicographic ordering}. 
\end{exam}
\noindent

Given a POVS $(V,V_+)$, denote
$$
  N_V:=\{x\in V|(\exists y\in V)(\forall n\in\nn)[-y\le nx\le y]\}
  \eqno(5)
$$
the set of all {\em infinitely small elements} of $V$. The wedge $V_+$ is almost Archi\-medean iff $N_V=\{0\}$. 
Since $V_+\cap (-V_+)\subseteq N_V$, we see that every almost Archi\-medean POVS is an OVS. Although, the wedge $V$ 
in a vector space $V\ne\{0\}$ is Archimedean but not almost Archi\-medean, any Archi\-medean OVS is clearly almost 
Archi\-medean. The converse is not true even if $\dim(V)=2$ (cf. also Example \ref{E2}).

\begin{exam}\label{E1}
Let $\Gamma$ be a set containing at least two elements, and $V$ be the space of all bounded real functions
on $\Gamma$, partially ordered by
$$
  f\le g \ \ \text{\it if} \ \ \text{\it either} \ f=g \ \  \text{\it or} \ \ \inf_{t\in\Gamma}[g(t)-f(t)]>0\,.
$$
The OVS $(V,\le)$ is almost Archi\-me\-dean but not Archi\-me\-dean. It can be shown that
$\inf_{n\ge 1}\frac{1}{n}f$ does not exist for any $0\ne f\in V_+$.  
\end{exam}

\begin{asser}\label{A2}
A wedge $K$ in a vector space $V$ is almost Archime\-dean if and only if it does not contain 
a straight line, say $x+ty$ with $x,y\in V$, $y\ne 0$. Furthermore, a POVS $(V,\le)$ 
is almost Archimedean if and only if
$$
  \bigcap_{n\ge 1}\Bigl[-\frac{1}{n}x,\frac{1}{n}x\Bigr]=\{0\}
$$
for every $x\in V_+$. 
\end{asser}
\noindent
The following assertion is immediate. 
\begin{asser}\label{A3}
Any subcone of an almost Archimedean cone is almost Archi\-medean.
\end{asser}

It is worth to notice, and it follows directly from Definitions \ref{D1}, \ref{D3}, that the intersection of 
any nonempty family of Archimedean (resp., almost Archimedean) wedges in a vector space is an Archimedean
(resp., almost Archimedean) wedge.

If $A$ is an order ideal in a POVS $(V,V_+)$ such that the POVS $(V/A,V_++A)$ is almost Archimedean, 
then $A$ is uniformly closed (cf. the part (a) of Proposition \ref{P2}). When $A=\{0\}$, the converse 
is true by the part (b) of Proposition \ref{P2}. Thus $V$ is almost Archimedean if and only if $\{0\}$ 
is uniformly closed. In general, the question, whether or not $V/A$ is Archimedean assuming an order 
ideal $A$ to be uniformly closed, is rather nontrivial (it has a positive answer due to Veksler \cite{V} 
in the vector lattice setting). An OVS $V$ is said to be a {\it vector lattice} (or a {\it Riesz space})
if every nonempty finite subset of $V$ has a least upper bound. It is well known that any almost 
Archimedean vector lattice is Archimedean. In Example \ref{E1}, $V=V/\{0\}$ is not Archimedean although $\{0\}$ 
is uniformly closed. Similarly to Definition \ref{D3}, we introduce almost Archimedean elements.

\begin{definition}\label{D4}
If $(V,V_+)$ is an POVS and $x\in V_+$, we say that $x$ is an {\rm almost Archimedean element} if,
whenever $y\in V$ with $x+ny\in V_+$ for all $n\in\zz$, then $y=0$.
\end{definition}

We shall use also the following definition.

\begin{definition}\label{D5}
Let $V=(V,V_+)$ be an OVS. Consider the category $\mathcal{C}_{Arch}(V)$ whose objects are pairs $\left\langle U,\phi\right\rangle$, 
where $U=(U,U_+)$ is an Archimedean OVS and $\phi:V\to U$ is a positive linear map $($that is: $\phi(V_+)\subseteq U_+$$)$, and morphisms 
$\left\langle U_1,\phi_1\right\rangle\to\left\langle U_2,\phi_2\right\rangle$ are positive linear maps $q_{12}: U_1\to U_2$ such that 
$q_{12}\circ \phi_1=\phi_2$. If $\mathcal{C}_{Arch}(V)$ possesses an initial object $\left\langle U_0,\phi_0\right\rangle$, then
$U_0$ is said to be an {\em Archimedeanization} of $V$.
\end{definition}

It should be noted that in Definitions \ref{D1}-\ref{D5} we use only the fact that $V$ is a {\em commutative ordered group},
where we replace the condition (1) in the definition of a wedge by
$$
  W\subseteq nW\ \ \ \ (\forall n\in\nn),
  \eqno(1')
$$ 
and positive linear maps in Definition \ref{D5} by group homomorphisms.

For further information on ordered vector spaces we refer to \cite{AT}.

\section{The Archimedeanization of an ordered vector space}

The following theorem is the main result of the present paper.

\begin{thm}\label{T1}
Any ordered vector space possesses a unique, up to an order isomorphism, Archime\-deanization.
\end{thm}

The proof is presented in the end of this section. First of all we introduce and discuss some
notions which may be considered as a generalization of the Archimedeanization method of 
Paulsen and Tomforde \cite{PT}. Given a POVS $(V,V_+)$, denote
$$
  D_V:=\{x\in V|(\exists \xi\in V_+)(\forall  n\in\nn)\ nx+\xi\in V_+\}\,,
$$
and let $I_V$ be the intersection of all uniformly closed order ideals in $(V,V_+)$.

It is straightforward to verify, that $D_V\cap(-D_V)=N_V$, where $N_V$ was defined in (5), 
$N_V$ is an order ideal in $(V,V_+)$, $D_V$ is a wedge, and $V_+\subseteq D_V$.
Consider the sets $V_++N_V=[V_+]_{N_V}$ and $D_V+N_V=[D_V]_{N_V}$ in the quotient space $V/N_V$.
Since $V_+$, $D_V$, and $N_V$ are wedges, then $V_++N_V$ and $D_V+N_V$ are wedges in $V/N_V$. 
However, in Example \ref{E0}, $D_V$ is not a cone. Clearly, 
$$
  (D_V+N_V)\cap -(D_V+N_V)=(D_V+N_V)\cap(-D_V+N_V)=
$$
$$  
  D_V\cap(-D_V)=N_V.
$$
So, $D_V+N_V$ is a cone in $V/N_V$, and $V_++N_V\subseteq D_V+N_V$ is a cone as well.
Let us collect some further elementary properties. 

\begin{prop}\label{P1}
Given an POVS $(V,V_+)$. Then

{\em (a)}\ $N_V$ is an order ideal in the POVS $(V,D_V)$$;$

{\em (b)}\ $N_V\subseteq I_V$$;$

{\em (c)}\ if $V$ possesses an order unit then $N_{V/A}=I_{V/A}$ for every ideal $A\subseteq V$$;$

{\em (d)}\ if $(V,V_+)$ is a Riesz space then $(V/N_V,D_V+N_V)$ is a Riesz space.
\end{prop}

{\bf Proof:}\ 
(a)\ Let $u,v\in N_V$, $w\in V$, with $u\le w\le v$ in $(V,D_V)$. Then 
$$
  -x\le nu\le x\,, \ \ -y\le nv\le y \ \ \ (\forall n\in\nn)  
  \eqno(6)
$$
for some $x,y\in V_+$, and 
$$
  0\le n(w-u)+\xi_1\,, \ \ 0\le n(v-w)+\xi_2 \ \ \ (\forall n\in\nn)  
  \eqno(7)
$$
for some $\xi_1,\xi_2\in V_+$. It follows from (6) and (7) that
$$
  0\le nw+x+\xi_1\,, \ \ 0\le -nw+y+\xi_2 \ \ \ (\forall n\in\nn).
$$
Then
$$
  -x+\xi_1\le nw\le y+\xi_2 \ \ \ (\forall n\in\nn)\,,
$$
and therefore
$$
  -(x+y+\xi_1+\xi_2)\le nw\le (x+y+\xi_1+\xi_2) \ \ \ (\forall n\in\nn)\,,
$$
which means that $w\in N_V$.

(b)\ Take an $f\in N_V$, then for some $g\in V$
$$
  -g\le nf\le g \ \ (\forall n\in\nn)
$$
or equivalently
$$
  -\frac{1}{n}g\le 0-f\le \frac{1}{n}g \ \ (\forall n\in\nn).
$$
Let $A\subseteq V$ be a uniformly closed order ideal. The sequence $a_n:=0\in A$ converges 
$g$-uniformly to $f$. Thus, $f\in A$. Since $A$ is arbitrary, then $f\in N_V$. 

(c)\ Denote by $e$ an order unit in $V$. Firstly, assume that $A=\{0\}$. Given $f\in N_V$, 
then for some $g\in V$
$$
  -g\le nf\le g \ \ (\forall n\in\nn)
$$
and, since $g\le ce$ for some $c\in\rr$, one gets 
$$
  -e\le nf\le e \ \ (\forall n\in\nn).
$$
We have to show that $N_V$ is uniformly closed. Suppose $f_n\stackrel{(u)}{\to}z$, then $f_n\stackrel{(e)}{\to}z$.
That is
$$
  -e \le k_n(f_n-z)\le e\ \ (\forall k\in\nn)
$$ 
for some $\nn\ni k_n\uparrow\infty$. But also 
$$
  -e \le k_nf_n\le e\ \ (\forall k\in\nn).
$$ 
Hence
$$
  -2e \le k_nz\le 2e\ \ (\forall k\in\nn).
$$ 
Since $k_n\uparrow\infty$, then $z\in N_V$.

Now, let $A\subseteq V$ be an arbitrary ideal. Since $e$ is an order unit in $V$, then $[e]$ is an order unit in $V/A$. 
So, $N_{V/A}=I_{V/A}$ as it was shown above.

(d)\ Denote by $u\vee v=\sup(u,v)$ and by $u\wedge v=\inf(u,v)$ the supremun and infimum of $u,v\in V$ in $(V,V_+)$. 
Take $u,v\in V$. Since $u\vee v -u, u\vee v -v \in V_+$ and $V_+\subseteq D_V$, then $u\vee v -u\in D_V$ and $u\vee v -v\in D_V$.

Let $w-u\,, w-v\in D_V$, then 
$$
  0\le n(w-u)+\xi_1\,, \ \ 0\le n(w-v)+\xi_2 \ \ \ (\forall n\in\nn)  
$$
for some $\xi_1,\xi_2\in V_+$. Denoting $\xi=\xi_1+\xi_2$, we obtain
$$
  0\le(n(w-u)+\xi)\wedge(n(w-v)+\xi)=n((w-u)\wedge(w-v))+\xi=
$$
$$
  n(w+(-u)\wedge(-v))+\xi=n(w-u\vee v)+\xi\ \ \ (\forall n\in\nn)\,, 
$$
which means that $w-u\vee v\in D_V$. As $N_V=D_V\cap(-D_V)$, the supremum of $[u],[v]$ exists, and is equal to $[u\vee v]$ in
the quotient OVS $(V/N_V,D_V+N_V)$. Since $u,v\in V$ are taken arbitrary, $(V/N_V,D_V+N_V)$ is a Riesz space.~$\blacksquare$

\begin{prop}\label{P2}
Given a POVS $(V,V_+)$. The following assertions hold true.

{\em (a)}\ Let $A\subseteq V$ be an order ideal. Then
$$
  V/A\ \ \text{\it is} \ \ \text{\rm almost}\ \ \text{\it Archimedean} \ \Rightarrow \ A \ \text{\it is} \ \
  \text{\it uniformly} \ \ \text{\it closed}\,.
$$

{\em (b)}\ \ $V \text{\it is} \ \ \text{\it almost}\ \ \text{\it Archimedean} \ \Leftrightarrow \ \ N_V=\{0\} \ \ \Leftrightarrow I_V=\{0\}$.

{\em (c)}\ If $(V,V_+)$ possesses an order unit then $(V/N_V,D_V+N_V)$ is an Archimedean OVS.

{\em (d)}\ If $(V/N_V, D_V+N_V)$ is Archi\-me\-dean then $(V/N_V, V_++N_V)$ is almost Archi\-me\-dean. In particular,
if $(V,V_+)$ possesses an order unit then $(V/N_V, V_++N_V)$ is almost Archimedean.
\end{prop}

{\bf Proof:}\ 
(a)\ Take a sequence $(y_n)\subseteq A$ such that $y_n\stackrel{(u)}{\to}x$ for some $u\in V_+$, and $x\in V$. We may assume that
$$
  n(y_n-x)\in [-u,u] \ \ \ (\forall n\in\nn).
$$ 
Thus
$$
  n[-x]=n[y_n-x]\in [-[u],[u]] \ \ \ (\forall n\in\nn).
$$ 
Since $V/A$ is almost Archimedean, $[-x]=[0]$, and hence $x\in A$.

(b)\ The implication ``$V$ is almost Archimedean $\Rightarrow\ \{0\}$ is uniformly closed'' follows from $(a)$.

The implication ``$\{0\}$ is uniformly closed $\Rightarrow\ I_V=\{0\}$'' is true due to the definition of $I_V$.  

The implication ``$I_V=\{0\}\ \Rightarrow\ N_V=\{0\}$'' is true by the part $(b)$ of Proposition \ref{P1}. 

The implication ``$N_V=\{0\}\ \Rightarrow\ V \ \text{\rm almost}\ \text{\rm Archimedean}$'' follows from Definition \ref{D3}. 

(c)\ It has been proved in \cite[Thm.2.35]{PT}.

(d)\ It follows from Assertion \ref{A3}, as $V_++N_V$ is a subcone of an almost Archimedean cone $D_V+N_V$ in $V/N_V$.~$\blacksquare$\\

The following example shows that, for the cone $D_V+N_V$ in $V/N_V$ to be Archi\-me\-dean, existence of an order unit in $(V,V_+)$
is not essential.

\begin{exam}\label{E2}
Let $V$ be the vector space $\mathcal{P}(\rr)$ of all real polynomials on $\rr$ with the positive cone
$$
  V_+:=\{p\in V|(\forall t\in\rr) p(t)>0\}\cup\{p\ |p(t)\equiv 0\}.
$$ 
The OVS $(V,V_+)$ has no order unit, $N_V=\{0\}$, and $D_V:=\{p\in V|(\forall t\in\rr) p(t)\ge0\}$. 
Take $x\equiv 1$, $y=-t^2$ in $V$. Then 
$$
   x-ny\in V_+\ \ \ \ (\forall n\in\nn)
$$
but $y\not\in -V_+$. Thus $(V,V_+)$ is not Archime\-dean, however $N_V=\{0\}$.
The OVS $(V/N_V,D_V+N_V)=(V,D_V)$ is clearly Archi\-medean.

Notice that if we take $V=\mathcal{P}[0,1]$ with the same ordering, then $V$ becomes a non-Archime\-dean but almost 
Archimedean OVS with an order unit and uniformly clo\-sed order ideal $N_V=\{0\}$. Moreover, it is well known that
this OVS satisfies the {\em Riesz decomposition property}, however is not a Riesz space.     
\end{exam}

Now, consider the following example due to T. Naka\-yama (see, for instance, \cite[p.436]{LZ}), in which
the cone $D_V+N_V$ in $V/N_V$ is not even almost Archi\-me\-dean. 

\begin{exam}[T. Nakayama]\label{E3}
Let
$$
  V=\{a=(a^1_k,a^2_k)_k|\ (a^1_k,a^2_k)\in(\rr^2,\le_{lex})\,, 
  a^1_k\ne 0   \ \text{\it for}\ \text{\it finitely}\ \text{\it many}\ k\}\,.
$$
Then $V$ is a Riesz space with respect to the pointwise ordering and operations. Denote
$$
  A=\{a\in V|\ a^1_k=0 \ \text{\it for}\ \text{\it all}\ k\,, a^2_k\ne 0 \ \text{\it for} \ \text{\it finitely}\ \text{\it many}\ k\}\,,
$$
$$
  \ \text{\it and}\ \ B:=\{a\in V|\ a^1_k=0 \ \text{\it for}\ \text{\it all}\ k\}\,.
$$
Clearly $A\subseteq B$, $A\ne B$, and $A\subseteq N_V$. Furthermore, $B=I_V$ and $A=N_V$. Indeed, 
take an $a\in B$, then 
$$
  -u \le n(a-a_n)\le u\ \ (\forall n\in\nn)
$$ 
for
$$
  a_n:=\{(0,a^2_1),(0,a^2_2),\ldots,(0,a^n_1),(0,0,),(0,0,),\ldots\}\in A 
$$ 
and $u:=(0,k|a^2_k|)_{k\in\nn}\in V$. So, $a_n\stackrel{(u)}{\to}a$ and $a\in I_V$, since $a_n\in N_V\subset I_V$. 
Therefore, $B\subseteq I_V$. Denote $M_k=\{a\in V|\ a^1_k=0\}$ for $k\in\nn$. 
Clearly, any $M_k$ is a uniformly closed ideal in $V$, and hence 
$$
  I_V\subseteq\bigcap_{k\in\nn}M_k=B\,.
$$  
So, $B=I_V$. In particular, $N_V\subseteq B$. But if $a\in N_V$ then
$$
  -b \le na\le b\in V\ \ \ (\forall n\in\nn)\,,
$$ 
and $b=(b^1_k,b^2_k)_{k\in\nn}$ with only finitely many $b^1_k\ne 0$. We may assume that if\ \ $b^1_k=0$ 
then also $b^2_k=0$. Hence $a=(a^1_k,a^2_k)_{k\in\nn}$ has all $a^1_k=0$ and only finitely many $a^2_k\ne 0$
Therefore $a\in A$ and, finally, $A=N_V$. Notice that 
$$
  N_V=\{a\in V|\ (\forall k) a^1_k=0 \ \text{\it and}\ a^2_k\ne 0 \ \text{\it for}\ \text{\it finitely}\ \text{\it many}\ k\}.
$$
$$
  D_V=\{a\in V|\ (\forall k) a^1_k\ge 0 \ \text{\it and}\ a^2_k<0 \ \text{\it for}\ \text{\it finitely}\ \text{\it many}\ k\}.
$$
Thus, $(V/N_V,D_V+N_V)$ is not necessarily almost Archimedean even if $V$ is a Riesz space. To see this, consider $c_n\in V$ with
$$
  \bigg[c_n\bigg]_k=
  \begin{cases}
  (0,-\frac{1}{2^k}) & 1 \le k \le n \\
  (0,0)      & k>n
  \end{cases}.
$$
Then $c_n\in N_V$, and $c_n$ converges uniformly to $d\in V$, where
$$
  d=\bigg((0,-\frac{1}{2^k})\bigg)_{k=1}^{\infty}\not\in N_V.
$$ 
Thus, $N_V$ is not uniformly closed in $(V,V_+)$. Moreover, $c_n$ converges uniformly to $d$ in $(V,D_V)$, 
however $d\not\in D_V$. In view of part $(2)$ of Proposition $\ref{P1}$, $(V/N_V,D_V+N_V)$ is a Riesz space. 
So, by Veksler's result {\em \cite{V}}, the Riesz space $(V/N_V,D_V+N_V)$ is not Archime\-dean and therefore  
not almost Archimedean.
\end{exam}

In connection with Example \ref{E2}, a question arises, whether or not any almost Archimedean 
OVS $(V,V_+)$ satisfies the property that $(V,D_V)$ is Archimedean. Fortunately, it has a positive
answer as the following proposition shows.

\begin{prop}\label{P3}
Let $(V,V_+)$ be an almost Archimedean OVS. Then the OVS $(V,D_V)$ is Archimedean.
\end{prop}

{\bf Proof:}\ 
Assume in contrary that the cone $D_V$ is not Archimedean. Then for some $x,y\in V$, $x\in D_V$,
$$
  (\forall r>0)\,\   rx+y\in D_V \ \ \text{\rm but}\ \ y\not\in D_V\,. 
$$
Since $x\in D_V$, then for some $v\in V_+$
$$
  (\forall \alpha>0)\, \ x+\alpha v\in V_+\,.
$$
Take
$$
  K=\spa (x,v,y)\cap V_+\,.
$$
Denote by $cl(C)$ the closure in Euclidean topology on $\spa(x,v,y)$ of $C\subseteq\spa(x,v,y)$. 
Then $x\in\cl(K)$ and $rx+y\in\cl(K)$ for all $r>0$ and hence $y\in \cl(K)\subseteq D_V$
which is impossible. Therefore, $D_V$ is Archimedean.~$\blacksquare$

\begin{coro}\label{C1}
Any almost Archimedean cone can be embedded into an Archime\-dean cone.
Moreover, given a POVS $(V,V_+)$, then the OVS $(V/N_V,D_V+N_V)$ is Archimedean if and 
only if the OVS $(V/N_V,V_++N_V)$ is almost Archimedean.
\end{coro}

{\bf Proof:}\ The first part is immediate by Proposition \ref{P3}. For the second part, denote $N=N_V$ and $D=D_V$. 

The necessity follows from Assertion \ref{A3}, since $V_++N$ is a subcone of the Archimedean cone $D+N$ in $V/N$.

For the sufficiency, assume that $(V/N,V_++N)$ is almost Archi\-medean. 
Then the cone $D_{V/N}$ is Archimedean in $V/N$, by Proposition \ref{P3}. Since 
$$
  D+N=\bigg[D\bigg]_N=\bigg[\{u\in V|(\exists\ \xi\in V_+)(\forall n\in\nn)\ nu+\xi\in V_+\}\bigg]_N\subseteq 
$$
$$
  \{[u]\in V/N|(\exists\ [\xi]\in V_++N)(\forall n\in\nn)[nu+\xi]\in V_++N\}=D_{V/N}, 
$$
the cone $D+N$ is almost Archimedean in $V/N_V=V/N$.~$\blacksquare$

Although, the first part Corollary \ref{C1} gives already an Archimedeanization of any almost Archimedean space $(V,V_+)$
(one only needs to replace $V_+$ by the intersection of all Archimedean cones containing $V_+$), our main result in this paper,
Theorem \ref{T1} do not require any restrictions on the cone $V_+$ and it covers also such OVS as $V$ in Example \ref{E3}
where $V/N$ still has nonzero infinitely small elements. 

{\bf Proof of Theorem \ref{T1}:}\ Let $V=(V,V_+)$ be an OVS. Denote $N_0=\{0\}$, 
$$
  N_1=N_V=\{x\in V|[x]_{N_0}=x\ \text{\rm is}\ \text{\rm infinitely}\ \text{\rm small}\ \text{\rm in}\ V/N_0=V\},
$$
$$
  N_{n+1}=\{x\in V|[x]_{N_n}\ \text{\rm is}\ \text{\rm infinitely}\ \text{\rm small}\ \text{\rm in}\ V/N_n\},
$$
and more generally
$$
  N_{\alpha}=N_{\alpha}(V)=\{x\in V|[x]_{\cup_{\beta<\alpha}N_{\beta}}\ \text{\rm is}\ \text{\rm infinitely}\ 
  \text{\rm small}\ \text{\rm in}\ V/_{\cup_{\beta<\alpha}N_{\beta}}\}
$$
for an arbitrary ordinal $\alpha>0$. It follows directly from the definition of $N_{\alpha}$ that
$$
  N_{\alpha_1}\subseteq N_{\alpha_2} \ \ \ (\forall \alpha_1\le\alpha_2). 
$$
Take the first ordinal, say $\lambda_V$, such that $N_{\lambda_V+1}=N_{\lambda_V}$. Then the OVS\\ 
$(V/N_{\lambda_V},[V_+]_{N_{\lambda_V}})$ has no nonzero infinitely small elements and hence  
is almost Archimedean. By Corollary \ref{C1}, the OVS 
$$
  (V/N_{\lambda_V},[D_{V/N_{\lambda_V}}]_{N_{\lambda_V}})=(V/N_{\lambda_V},D_{V/N_{\lambda_V}})
$$ 
is Archimedean. Denote by $p_V$ the quotient map $V\to V/N_{\lambda_V}$, then $p_V$ is positive and linear. 
For any other pair $\left\langle U,\phi \right\rangle$, where $U=(U,U_+)$ is an Archimedean OVS and 
$\phi:V\to U$ is a positive linear map, we have that $\phi(x)\in N_U$ for every $x\in N_{\alpha}$,
where $\alpha$ is an arbitrary ordinal. Since every Archimedean OVS is almost Archimedean, $N_U=\{0\}$ and therefore
$N_{\lambda_V}\subseteq\ker(\phi)$. So, the map $\tilde{\phi}:V/N_{\lambda_V}\to U$ is well defined 
by $\tilde{\phi}([x]_{N_{\lambda_V}})=\phi(x)$ and satisfies $\tilde{\phi}\circ p_V=\phi$. 
Moreover, if $[x]_{N_{\lambda_V}}\in D_{V/N_{\lambda_V}}$, then  
$$
  n[x]_{N_{\lambda_V}}+[\xi]_{N_{\lambda_V}}\in [V_+]_{N_{\lambda_V}} \ \ \ (\forall n\in\nn)  
$$
for some $[\xi]_{N_{\lambda_V}}\in [V_+]_{N_{\lambda_V}}$. Applying $\phi$ gives
$$
  n\phi(x)+\phi(\xi)\in U_+ \ \ \ (\forall n\in\nn),  
$$
and, since $(U,U_+)$ is Archimedean, it follows that $\phi(x)\in U_+$. Thus
$$
  \tilde{\phi}([x]_{N_{\lambda_V}})=\phi(x)\in U_+\,,
$$
and $\tilde{\phi}$ is a positive linear map. To show that $\tilde{\phi}$ is unique, take any $\psi:V/N_{\lambda_V}\to U$, that 
satisfies $\psi\circ p_V=\phi$. Then 
$$
  \psi([y]_{N_{\lambda_V}})=\psi(p_V(y))=\phi(y)=\tilde{\phi}(p_V(y))=\tilde{\phi}([y]_{N_{\lambda_V}})\ \ \ \ (\forall y\in V)\,,
$$
and hence $\psi=\tilde{\phi}$. Thus, $\left\langle(V/N_{\lambda_V},D_{V/N_{\lambda_V}}),p_V\right\rangle$ 
is an initial object of the category $\mathcal{C}_{Arch}(V)$. Hence, the OVS $(V/N_{\lambda_V},D_{V/N_{\lambda_V}})$ 
is an {\em Archimedeanization} of the OVS $(V,V_+)$. 
~$\blacksquare$

In connection with the proof of Theorem \ref{T1}, the following question arises naturally.
Whether or not for any ordinal $\alpha$ there is an OVS $(V,V_+)$ for which $\alpha$
is the first ordinal such that $N_{\alpha+1}(V)=N_{\alpha}(V)$?

\end{document}